\documentclass[letterpaper,final,twoside,reqno]{amsart}
\usepackage{xspace}
\usepackage[all]{xy}
\usepackage{hyperref}

\newcommand{\baseRing}[1]{\ensuremath{\mathbb{#1}}}
\newcommand{\Z}{\baseRing{Z}}
\newcommand{\Q}{\baseRing{Q}}
\newcommand{\N}{\baseRing{N}}
\newcommand{\R}{\baseRing{R}}
\newcommand{\C}{\baseRing{C}}
\newcommand{\CP}{\baseRing{P}}
\newcommand{\jdef}[1]{\emph{#1}}
\newcommand{\pd}[2]{{\frac{\partial {#1}}{\partial {#2}}}}

\newcommand{\conj}{\overline}
\newcommand{\DD}{\ensuremath{{\mathcal D}}\xspace}

\theoremstyle{plain}
\newtheorem{theorem}{Theorem}[section]
\newtheorem{corollary}[theorem]{Corollary}
\newtheorem{proposition}[theorem]{Proposition}
\newtheorem{lemma}[theorem]{Lemma}

\theoremstyle{definition}
\newtheorem{definition}[theorem]{Definition}
\newtheorem{remark}[theorem]{Remark}
\newtheorem{example}[theorem]{Example}

\numberwithin{equation}{section}

\DeclareMathOperator{\symm}{Symm}
\DeclareMathOperator{\vol}{vol}
\DeclareMathOperator{\lcm}{lcm}
\DeclareMathOperator{\End}{End}

\newcommand{\IVHS}{infinitesimal variation of Hodge structure\xspace}
\newcommand{\IVHSs}{infinitesimal variations of Hodge structure\xspace}

\CompileMatrices

\begin{document}

\title[Non-genericity of IVHS]{Non-genericity of infinitesimal
  variations of Hodge structures arising in some geometric contexts}
\author{Emmanuel Allaud}
\author{Javier Fernandez}

\address{Department of Mathematics\\ University of Utah\\ Salt Lake
  City\\ UT 84112--0090\\USA}
\email{allaud@math.utah.edu}

\address{Instituto Balseiro\\ Universidad Nacional de
  Cuyo -- C.N.E.A.\\Bariloche\\R8402AGP\\Rep\'ublica Argentina}
\email{jfernand@ib.edu.ar}


\bibliographystyle{amsplain}


\begin{abstract}
  We prove that the infinitesimal variations of Hodge structure
  arising in a number of geometric situations are non-generic. In
  particular, we consider the case of generic hypersurfaces in
  complete smooth projective toric varieties, generic hypersurfaces in
  weighted projective spaces and generic complete intersections in
  projective space and show that, for sufficiently high degrees, the
  corresponding infinitesimal variations are non-generic.
\end{abstract}

\maketitle



\section{Introduction}
\label{sec:intro}

A variation of Hodge structure can be described, using the language of
exterior differential systems, as an integral manifold of Griffiths'
differential system over the period domain. An important problem in
Hodge theory is the study of the geometric locus, that is, the locus
of those variations of Hodge structure that arise from the cohomology
of a family of polarized projective varieties. 

An infinitesimal version of this problem consists of describing the
infinitesimal variations of Hodge structure ---the integral elements
of Griffiths' system--- that arise from geometric variations.
In~\cite{ar:allaud-hypersurfaces_dn} the first author has shown that
the infinitesimal variations arising from deformations of
hypersurfaces of sufficiently high degree in projective space are
non-generic in the space of infinitesimal variations. The purpose of
the present paper is to show that this property holds in a variety of
geometric situations, namely, generic hypersurfaces of complete smooth
toric varieties (Theorem~\ref{th:smooth_main}), generic hypersurfaces
in weighted projective space (Theorem~\ref{th:weighted_main}) and
generic complete intersections in projective space
(Theorem~\ref{th:ci_main}).  In all cases, requirements of
sufficiently high degree apply.

Considering all these results we begin to see a general principle that
infinitesimal variations of geometric origin (eventually satisfying
some condition analogous to high degree) are non-generic.

The main tools used in this paper are the appropriate residue theories
for simplicial toric varieties and for complete intersections in
projective space, as well as infinitesimal Torelli theorems, dualities
and Macaulay's theorem.

The plan for the paper is as follows: in Section~\ref{sec:ideas} we
review some results from~\cite{ar:allaud-hypersurfaces_dn} and
describe an approach to proving the non-genericity of families of
infinitesimal variations. In Section~\ref{sec:toric} we study the
infinitesimal variations associated to generic hypersurfaces in
complete simplicial projective toric varieties and reduce the proof of
non-genericity results to a numerical condition strongly related to
infinitesimal Torelli theorems
(Theorem~\ref{th:non_genericity_toric}). One such Torelli theorem, due
to M. Green, allows us to conclude the non-genericity in the case of
smooth ambient spaces.  In Section~\ref{sec:hypersurfaces_in_wpn} we
specialize the toric analysis to the case where the ambient space is a
weighted projective space where the numerical condition now follows
from a result of L. Tu. Finally, in Section~\ref{sec:ci_case}, the
general approach is specialized to the case of complete intersections.

Finally, we wish to thank D. Cox, E. Cattani, A. Dickenstein and J. M.
Landsberg for many helpful discussions.
  


\section{The projection of the integral elements of the Griffiths'
  system to a grassmannian, the symmetrizers correspondence}
\label{sec:ideas}

Let us fix some notation (see \cite{bo:griffiths-topics} and
\cite{ar:Mayer-coupled} for more details). Recall that a \jdef{(real)
  Hodge structure of weight $k$} on the real vector space $H$ with
Hodge numbers $h^{k,0},\ldots,h^{0,k}$, consists of a grading
$H\otimes \C = H^{k,0} \oplus \cdots \oplus H^{0,k}$ such that $\dim
H^{j,k-j} = h^{j,k-j}$ and $\conj{H^{j,k-j}} = H^{k-j,j}$. We will
usually denote $H$ and its complexification by $H$.  We define $H^q:=
\hom(H^{k-q,q},H^{k-q-1,q+1})$, and consider $\oplus_{0 \leq q \leq
  k-1} H^q$ as a subset of $\hom(H,H)$. We also define the following
maps for $0 \leq q \leq k-1$:
\begin{align*}
  p_q : \oplus_{0 \leq r \leq k-1} H^r &\longrightarrow H^q \\
  \alpha &\longmapsto \alpha_{\mid H^{k-q,q}}
\end{align*}

\begin{remark}
  $p_q$ is the natural projection of $\oplus_{0 \leq a \leq k-1} H^a$
  onto $H^q$.
\end{remark}

The periods space (i.e. the set of all \emph{polarized} Hodge
structures with fixed Hodge numbers and polarization $Q$,
see~\cite{bo:griffiths-topics}) is the homogeneous variety $\DD \simeq
G/P$ (with $G=SO(H,Q)$ and $P$ a parabolic subgroup). We note
$\mathfrak g$, $\mathfrak p$ the Lie algebras of $G$ and $P$.  Then
$\mathfrak g$ is given by:
\begin{equation}\label{eq_alg_lie}
  \mathfrak g \ =\ \left\{ X \in \End(H) \mid Q(X v, w) + Q(v, X w) = 0,\,
    \forall v,w \in H \right \}
\end{equation}
Fixing a reference structure $H_0 := \{H^{k-q,q}_0\} \in \DD$, we can also
consider the following subspaces of $End(H)$:
\begin{equation}\label{sh_sur_end_H}
  \End(H)^{p,-p}\ =\ \left\{ X \in \End(H) \mid \forall r+s=k,X(H^{r,s}_0)
  \subset H^{r+p,s-p}_0 \right\} 
\end{equation}
and then define $\mathfrak g^{p,-p} := \mathfrak g \cap \End(H)^{p,-p}$.
We also note
\begin{equation*}
  \mathfrak g^0 \ :=\  g^{0,0},\quad \mathfrak g^- \ :=\ \bigoplus_{\substack p<0}
  \mathfrak g^{p,-p},\quad \mathfrak g^+ \ :=\ \bigoplus_{\substack p>0} \mathfrak
  g^{p,-p}
\end{equation*}
So that we have $\mathfrak g=\mathfrak g^- \oplus \mathfrak g^0
\oplus \mathfrak g^+$. Moreover we have that $\mathfrak p=\mathfrak
g^0 \oplus \mathfrak g^+$ and, as $\DD \simeq G/P$, we conclude that
\begin{equation*}
  T_{H_0} \DD \ \simeq\  \mathfrak g^-.  
\end{equation*}
Let $V:=\mathfrak g^{-1,1} \subset \oplus_{0 \leq a \leq k-1} H^a$,
and let $V_d$ denote the set of \IVHSs\ of dimension $d$, then $V_d
\subset G(d,V)$. In fact, $V_d$ is the algebraic subvariety of $G(d,V)$
of all abelian $d$-subalgebras of $V$.  Moreover we have the
following:
\begin{lemma}
  For all $i \in \{0,\dots,k\}$, $p_i : G(d,V) \rightarrow G(d,H^i)$
  is a rational map.
\end{lemma}

We now recall the definition of Symmetrizer due to R. Donagi:
\begin{definition}
  Let $\psi:E \times F \rightarrow G$ be a bilinear map. We define
  \begin{equation*}
    \symm \psi \ :=\ \{ q \in \hom(E,F) \mid \forall\, \alpha,\alpha^\prime
    \in E, \psi(\alpha,q(\alpha^\prime))=\psi(\alpha^\prime,q(\alpha))
    \}     
  \end{equation*}
\end{definition}

The link between the \IVHSs\ and the symmetrizers is given by the
following proposition (\cite{ar:allaud-hypersurfaces_dn}, Proposition
3.2):
\begin{proposition}\label{proposition:fibre-proj1}
  For any $E^0 \in G(d,H^0)$, we define the following bilinear map:
  \begin{align*}
    \phi_{E^0} : E^0 \times H^1 &\longrightarrow \hom(H^{k,0},H^{k-2,2}) \\
    (\alpha,\beta) &\longmapsto \beta \circ \alpha
  \end{align*}
  then we have for $k\geq 3$
  \begin{equation*}
    p_1(p_0^{-1}(E^0) \cap V_d) \ \simeq\ \symm \phi_{E^0}    
  \end{equation*}
\end{proposition}

We need one more proposition to explain the principle of
``non-genericity'' for \IVHSs.
\begin{proposition} \label{proposition:symm_rank}
  Let $G^0$, $G^1$, $G^2$ be three $\C$-vector spaces with $\dim
  G^0>1$. For all $E^0 \subset \hom(G^0,G^1)$, we consider the
  following bilinear map :
  \begin{align*}
    \phi_{E^0} : E^0 \times \hom(G^1,G^2) &\longrightarrow \hom(G^0,G^2) \\
    (\alpha,\beta) &\longmapsto \beta \circ \alpha
  \end{align*}
  Let $d$ be a positive integer, $d \leq \dim \hom(G^0,G^1)$; if we
  have the following inequality
  \begin{equation} \label{ineq-symm-trivial}
    d \geq 3\left(\left[ \frac{\dim G^1-1}{\dim G^0} \right]+1\right)
  \end{equation}
  then any generic $E^0 \in G(d,\hom(G^0,G^1))$ satisfies
  \begin{equation*}
    \symm(\phi_{E^0})\ =\ \{0\}    
  \end{equation*}
\end{proposition}

\begin{proof}
  (\cite{ar:allaud-hypersurfaces_dn}, Proposition 3.3) The space of
  symmetrizers $\symm(\phi_{E^0})$ is defined by a set of linear
  equations depending on $E^0 \in G(d,\hom(G^0,G^1))$; then the
  condition $\symm(\phi_{E^0})=\{0\}$ is open in $G(d,\hom(G^0,G^1))$
  and the proof amounts to build a $\C$-vector space in
  $\hom(G^0,G^1)$ for which the rank of the linear equations defining
  the symmetrizers is maximal (the dimensions are such that maximal
  rank for those equations implies triviality of their zero set).
\end{proof}

In our context we will have $G^i:=H^{k-i,i}$ ($0 \leq i \leq 2$), and
$E^0:=p_0(E)$ where $E$ will be an \IVHS. In this setting the
inequality~\eqref{ineq-symm-trivial} now reads:
\begin{equation}\label{eq:ineq-principle}
  \dim p_0(E) \ \geq\  3\left(\left[ \frac{h^{k-1,1}-1}{h^{k,0}}
  \right]+1\right) 
\end{equation}

Now we can state the main theorem that we will use in our geometric
applications.

\begin{theorem} \label{th:non-genericity-principle}
  Let $E$ be an \IVHS\ such that $\dim p_0(E)=\dim E$, $p_1(E) \neq
  \{0\}$ and $E$ satisifies inequality~\eqref{eq:ineq-principle} then $E$
  is non-generic, that is, $E$ must lie in a proper subvariety of
  $V_{\dim E}$.
\end{theorem} 

\begin{proof}
  Define $E^0:=p_0(E)$ and $d:=\dim E$. Using Proposition
  \ref{proposition:fibre-proj1} we have that $p_1(E) \in
  p_1(p_0^{-1}(E^0) \cap V_d) \simeq \symm \phi_{E^0}$, and as we made
  the assumption that $p_1(E) \neq \{0\}$, this implies that $\symm
  \phi_{E^0} \neq \{0\}$. But as $E$ satisfies~\eqref{eq:ineq-principle},
  the generic element $X \in G(d,H^0)$ satisfies $\symm \phi_X=\{0\}$,
  then $p_0(E)$ must lie in a proper subvariety $S$ of $G(\dim p_0(E),
  H^0)$; therefore $E$ lies in $p_0^{-1}(S) \cap V_{\dim E}$ which is
  a proper subvariety of $V_d$ (because $p_0$ is regular and surjective).
\end{proof}

The following are applications of
Theorem~\ref{th:non-genericity-principle}; they are the main results
of this paper.

\begin{theorem}
  Let $P$ be a complete smooth projective toric variety of dimension
  $n\geq 4$ and $\beta\in A_{n-1}(P)$ the degree of an ample Cartier
  divisor. Let $E$ be the infinitesimal variation of Hodge structure
  associated to a generic $f \in S_{t\beta}$, where $S$ is the
  homogeneous coordinate ring of $P$. Then, there exists $t_0\in \N$
  such that for all $t\geq t_0$, $E$ is non-generic.
\end{theorem}

\begin{theorem}
  Let $P=\CP^n(q_0,\ldots,q_n)$ be the weighted projective space with
  $n\geq 4$, $q_0 = \gcd(q_1,\ldots, q_n) = 1$. Then, if
  $\lcm(q_1,\ldots,q_n)|d$ and $\lcm(q_1,\ldots,q_n)|\sum_{j=0}^n
  q_j$, there is a degree $d_0$ such that for all $d\geq d_0$ the
  infinitesimal variation of Hodge structure associated to a generic
  hypersurface of degree $d$ is non-generic.
\end{theorem}

\begin{theorem}
  Let $X\subset \CP^n$ be a smooth complete intersection defined by
  the intersection of $c$ hypersurfaces of degrees $d_1,\ldots,d_c$.
  If $\dim X \geq 3$ there exists $K\in\N$ such that if the degree of
  the canonical bundle $d(X) = \sum_{a=1}^c d_a -(n+1)\geq K$, then
  the \IVHS associated to the deformation of $X$, $E$, is non-generic
  in $V_{\dim E}$. 
\end{theorem}

\begin{remark}
  These theorems state asymptotic results.  Nevertheless, for concrete
  examples the conditions on the degrees of the hypersurfaces that
  ensure non-genericity can be stated and tested explicitly. For the
  last Theorem, effective bounds are also included.
\end{remark}


\section{Hypersurfaces in toric varieties}
\label{sec:toric}

In this section we describe the infinitesimal variation of Hodge
structure associated to the family of hypersurfaces of fixed degree in
a complete simplicial projective toric variety $P$. Using this
description the non-genericity of the infinitesimal variation for
hypersurfaces of high degree can be established modulo some algebraic
conditions. Finally, we show that the algebraic conditions hold, for
instance, when the hypersurfaces and $P$ are smooth and, also, when
$P$ is a weighted projective space.

Let $P$ be the $n$-dimensional complete simplicial toric variety
defined by the fan $\Sigma\subset N$, where $N$ is a lattice of rank
$n$ in $\R^n$. $P$ can be described as a geometric quotient as
follows~\cite{ar:cox-homogeneous}. Let $S:=\C[z_1,\ldots,z_r]$ be the
polynomial ring with variables corresponding to the integer generators
$n_j$ of the $1$-dimensional cones of $\Sigma$. Let $Z:= \cap_{\sigma
  \subset \Sigma} \{z\in \C^r : \hat{z}_\sigma = 0\}$, where for each
cone $\sigma \subset \Sigma$, $\hat{z}_\sigma := \prod_{n_j \notin
  \sigma} z_j$. Then $P$ is the geometric quotient of $U:=\C^r
\setminus Z$ by the algebraic group $G:=\hom(A_{n-1}(P),\C^*)$, where
$A_{n-1}(P)$ is the Chow group of $P$.

If $M$ is the lattice dual to $N$, a grading on $S$ can be induced
from the exact sequence $0\rightarrow M \rightarrow \Z^r \rightarrow
A_{n-1}(P) \rightarrow 0$, where the second arrow is $m\mapsto
(\langle m,n_1\rangle, \ldots, \langle m,n_r\rangle)$ and the third is
$(a_1,\ldots,a_r)\mapsto a_1 [D_1] + \cdots + a_r [D_r]$, where $D_j$
is the divisor associated to the generator $n_j$. With this notation,
$\deg(z_1^{a_1} \cdots z_r^{a_r}) := a_1 [D_1] + \cdots + a_r [D_r]
\in A_{n-1}(P)$. This grading coincides with the one induced by the
action of $G$ on $\C^r$. The graded piece of $S$ of degree $\beta \in
A_{n-1}(P)$ is denoted by $S_\beta$.

Given $f\in S_\beta$, let $V(f):=\{z \in \C^r : f(z)=0\}$; $V(f)\cap
U$ is $G$ stable and hence descends to a hypersurface $X_f\subset P$.
Since $\Sigma$ is simplicial, $P$ is a $V$-manifold. Moreover, by
Proposition 4.15 in~\cite{ar:batyrev_cox-hypersurfaces}, if $\beta$ is
Cartier and ample, for generic $f\in S_\beta$, $X_f$ is a
$V$-submanifold of $P$.  Let $B_\beta\subset S_\beta$ be the Zariski
open of those $f$ for which $X_f$ is a $V$-submanifold of $P$ and let
$\pi:{\mathcal X}\rightarrow B_\beta$ be the family of hypersurfaces
${\mathcal X}_f = X_f \subset P$. If we assume that $P$ is a
projective variety, classical results of P. Deligne show that
$R^{n-1}\pi_*\Q$ defines a graded-polarized variation of mixed Hodge
structure. On the other hand, by~\cite[\S
1]{ar:steenbrink-vanishing_cohomology}, it is known that since $X_f$
is a $V$-manifold, $H^{n-1}(X_f,\Q)$ carries a pure Hodge structure of
weight $n-1$. Furthermore, if $\omega:=c_1({\mathcal O}_P(1))$ and
$L_\omega$ is the operator induced by left multiplication by $\omega$,
\begin{equation*}
  H^{n-1}(X_f)_{p} \ :=\ \ker(L_\omega: H^{n-1}(X_f,\Q) \rightarrow
  H^{n+1}(X_f,\Q)),
\end{equation*}
the \jdef{primitive cohomology}, is a sub Hodge structure that is
polarized. In what follows we will consider a sub Hodge structure of
$H^{n-1}(X_f)_{p}$, namely, the \jdef{vanishing cohomology} (also
known as \jdef{variable cohomology}) defined as
\begin{equation*}
  H^{n-1}(X_f)_{v}\ :=\ \ker(j_* : H^{n-1}(X_f,\Q)\rightarrow
  H^{n+1}(X_f,\Q)),
\end{equation*}
where $j_*$ is the Gysin morphism associated to the inclusion
$j:X_f\rightarrow P$. Notice that
in~\cite{ar:batyrev_cox-hypersurfaces} this vanishing cohomology is
called primitive cohomology. Since $L_\omega = j^*\circ j_*$ it is
clear that $H^{n-1}(X_f)_{v}$ is a sub Hodge structure of
$H^{n-1}(X_f)_{p}$, that is polarized, hence, $H^{n-1}(X_f)_{v}$ is a
polarized sub Hodge structure of $H^{n-1}(X_f)_{p}$ of weight $n-1$.
To sum up, we have:

\begin{proposition}\label{prop:pvhs_van}
  For $\beta\in A_{n-1}(P)$ Cartier and ample, the vanishing
  cohomology of the family of hypersurfaces $\pi:{\mathcal
    X}\rightarrow B_\beta$ defines the polarized variation of Hodge
  structure of weight $n-1$, whose fibers are $H^{n-1}(X_f,\Q)_{v}$.
\end{proposition}

If $\DD$ is the classifying space for the Hodge structures described
in Proposition~\ref{prop:pvhs_van} and $\Gamma$ is the monodromy
group, the previous result defines a \jdef{period mapping}
$\Phi:B_\beta\rightarrow \DD/\Gamma$. In fact, upon restriction to a
Zariski open $B_\beta'\subset B_\beta$, $\Phi$ factors through the
\jdef{generic coarse moduli space for hypersurfaces of $P$ with degree
  $\beta$}, $M_\beta:= B_\beta'/\widetilde{Aut_\beta}$, as defined in~\cite[\S
13]{ar:batyrev_cox-hypersurfaces}:
\begin{equation}\label{diag:factorization_through_moduli}
  \xymatrix{{B'_\beta}\ar@{^{(}->}[r] \ar[d]_{\kappa} 
    & {B_\beta} \ar[r]^{\Phi}
    & {\DD/\Gamma}\\
    {M_\beta} \ar[rru]_\Psi & & }
\end{equation}
As noted in~\cite{ar:batyrev_cox-hypersurfaces}, for $B'_\beta$
sufficiently small, $\kappa$ is a smooth map. Taking differentials at
the point $f\in B'_\beta$ we have $d\Phi_f = d\Psi_{[f]} \,
d\kappa_f$, with $d\kappa_f$ onto. Notice that the horizontality of
$\Phi$ implies that of $\Psi$. In what follows we will consider the
infinitesimal variation of Hodge structure $E_f$ associated to $\Psi$,
$d\Psi_{[f]}:T_{[f]}M_\beta \rightarrow T_{\Psi([f])}\DD$. In fact, we
will be only interested in the image of $d\Phi_f$, that we will
continue to call $E_f$. Since $d\kappa_f$ is onto, $E_f$ is also the
image of $d\Psi_{[f]}$.

The following result shows some relations between the cohomology of a
hypersurface $X_f$ in a toric variety and graded parts of the Jacobian
ring of $f$.

\begin{theorem}\label{th:toric_hypersurfaces_theory}
  Let $P$ be a complete simplicial projective toric variety of
  dimension $n\geq 4$, and $X_f\subset P$ a quasi-smooth ample
  hypersurface defined by $f\in S_\beta$ that defines a Cartier
  divisor. Let $R(f):= S / \langle \pd{f}{z_1}, \ldots, \pd{f}{z_r}
  \rangle$ be the Jacobian ring of $f$, with the grading inherited
  from $S$. Then,
  \begin{enumerate}
  \item \label{it:compatibility_2} $H^0(P,{\mathcal O_P(X_f)})\simeq
    S_{\beta}$. The following diagram is commutative up to a
    multiplicative constant.
    \begin{equation}
      \label{diag:compatibility_toric_2}
      \xymatrix{{H^0(P,{\mathcal O}_P(X_f)) \otimes H^0(P,\Omega^n_P(p X_f))} 
        \ar[r] \ar@{->>}[d] & 
        {H^0(P,\Omega^n_P((p+1) X_f))} \ar@{->>}[d] \\
      {T_fB_\beta \otimes H^{n-p,p-1}(X_f)_v} \ar[r] & {H^{n-p-1,p}(X)_v}}
    \end{equation}
    where the top arrow is induced by the product of sections and the
    bottom by the infinitesimal variation of Hodge structure, and the
    vertical arrows are the surjective maps 
    \begin{equation*}
      \bar{\alpha}_p:H^0(P,\Omega^n_P(pX_f))\rightarrow H^{n-p,p-1}(X_f)_v 
    \end{equation*}
    induced by the residue.
  \item \label{it:compatibility} The diagram
    \begin{equation}\label{diag:compatibility_toric}
      \xymatrix{{T_{[f]}M_\beta} \ar[r]^-{d\Psi_{[f]}} \ar[d]_{\simeq} &
        {\hom(H^{n-1,0}(X_f)_{v}, H^{n-2,1}(X_f)_{v})} \ar[d]^{\simeq} \\
        {R_\beta} \ar[r]^-{\times} & {\hom(R(f)_{\beta -\beta_0},
        R(f)_{2\beta -\beta_0})} }
    \end{equation}
    is commutative up to a multiplicative constant. Here $\beta_0 :=
    \deg(z_1\ldots z_r)\in A_{n-1}(P)$.
  \end{enumerate}
\end{theorem}

\begin{proof}
  The first assertion of~\eqref{it:compatibility_2} is Lemma 4.11
  in~\cite{ar:batyrev_cox-hypersurfaces}. Theorem 6.13
  from~\cite{bo:Voisin_Hodge_II} (modified for toric varieties using
  the residue theory developed in~\cite{ar:batyrev_cox-hypersurfaces})
  gives that the diagram
  \begin{equation}\label{diag:6.13}
    \xymatrix{{H^0(P,\Omega^n_P(p X_f))} \ar@{->>}[d] \ar[r] & 
      {\hom(H^0(P,{\mathcal O}_P(X_f)),H^0(P,\Omega^n_P((p+1) X_f)))} 
      \ar@{->>}[d] \\ {H^{n-p,p-1}(X_f)_v} \ar[r] & 
      {\hom(T_f B_\beta, H^{n-p-1,p}(X_f)_v)}}
  \end{equation}
  is commutative up to a multiplicative constant. The vertical arrows
  are induced by the surjective residue maps $\bar{\alpha}_p$ and the
  isomorphism $H^0(P,{\mathcal O_P(X_f)})\simeq S_{\beta} \simeq T_f B_\beta$.
  From~\eqref{diag:6.13} the rest of part~\eqref{it:compatibility_2}
  follows.

  To prove part~\eqref{it:compatibility}, we start by showing that the
  vertical arrows are isomorphisms. Proposition 13.7
  in~\cite{ar:batyrev_cox-hypersurfaces} does it for the left arrow,
  while Theorem 10.13 in~\cite{ar:batyrev_cox-hypersurfaces} implies
  the result for the right one.

  Next we establish the commutativity of
  diagram~\eqref{diag:compatibility_toric}. Starting
  from~\eqref{diag:6.13} and using Lemma 4.11 and Theorem 9.7
  of~\cite{ar:batyrev_cox-hypersurfaces}, we have that
  \begin{equation*}
    \xymatrix{ {S_{p\beta - \beta_0}} \ar[d]_{\simeq} \ar[r] & 
      {\hom(S_\beta, S_{(p+1)\beta - \beta_0})} \ar[d]^{\simeq} \\      
      {H^0(P,\Omega^n_P(p X_f))} \ar@{->>}[d] \ar[r] & 
      {\hom(H^0(P,{\mathcal O}_P(X_f)),H^0(P,\Omega^n_P((p+1) X_f)))} 
      \ar@{->>}[d] \\ {H^{n-p,p-1}(X_f)_v} \ar[r] & 
      {\hom(T_f B_\beta, H^{n-p-1,p}(X_f)_v)}}
  \end{equation*}
  commutes up to a multiplicative constant. From here it is easy to
  see that the same holds for
  \begin{equation*}
    \xymatrix{{R_{p\beta - \beta_0}} \ar[d]_{\simeq} \ar[r] & 
    {\hom(S_\beta, R_{(p+1)\beta - \beta_0})} \ar[d]^{\simeq} \\ 
    {H^{n-p,p-1}(X_f)_v} \ar[r] & 
    {\hom(T_f B_\beta, H^{n-p-1,p}(X_f)_v)}}
  \end{equation*}
  and, eventually, for
  \begin{equation*}
    \xymatrix{{S_\beta} \ar[r] \ar[d]_{\simeq} & {\hom(R_{p\beta-\beta_0}, 
        R_{(p+1)\beta-\beta_0})} \ar[d]^{\simeq} \\ {T_f B_\beta} \ar[r] & 
        {\hom(H^{n-p,p-1}(X_f)_v, H^{n-p-1,p}(X_f)_v)}}
  \end{equation*}
  and, then, \eqref{diag:compatibility_toric} follows by noticing that
  the horizontal arrows factor through $R_\beta$ and $T_{[f]}M$
  respectively.
\end{proof}

V. Batyrev gives a similar description of the differential of the
period mapping in his work on affine hypersurfaces (see Proposition
11.8 in~\cite{ar:Batyrev_affine_VMHS}).

To state the following result we will write $O(t^k)$ to denote a
function ---a polynomial in this case--- with $\frac{O(t^k)}{t^k}$
bounded as $t\rightarrow \infty$. Also, if $M\simeq \Z^n \subset \R^n$
is a lattice, $\vol(\pi)$ is the volume of the lattice polytope
$\pi\subset \R^n$, normalized so that the unit $n$-cube of the lattice
$M$ is $1$.

\begin{lemma}\label{le:dimensions_toric}
  Let $D$ be an ample Cartier divisor in the $n$-dimensional complete
  simplicial toric variety $P$ and $\Delta\subset \R^n$ its associated
  polytope. If $X$ is a generic ample hypersurface of degree $[tD]\in
  A_{n-1}(P)$ for $t\in \N$, then $\vol(\Delta) > 0$ and
  \begin{enumerate}
  \item $h^{n-1,0}_t := \dim H^{n-1,0}(X) = \vol(\Delta) t^n +
    O(t^{n-1})$,
  \item $h^{n-2,1}_t := \dim H^{n-2,1}(X) = (2^n-(n+1))\vol(\Delta)
    t^n + O(t^{n-1})$ and
  \item $\mu_t:=\dim M_{[tD]} \geq \vol(\Delta) t^n + O(t^{n-1})$.
  \end{enumerate}
\end{lemma}

\begin{proof}
  Let $\Delta_t := t\Delta$ be the polytope in $M\otimes\R \simeq
  \R^n$ defined by $tD$. By (5.5)
  in~\cite{ar:danilov_khovanski-newton} we have
  \begin{equation}\label{eq:h_computation}
    \begin{split}
      h^{n-1,0}_t \ &:=\  \dim H^{n-1,0}(X) \ =\  l^*(\Delta_t) \\
      h^{n-2,1}_t \ &:=\ \dim H^{n-2,1}(X) \ =\  l^*(2\Delta_t) - (n+1)
      l^*(\Delta_t) - \sum_{\Gamma\in {\mathcal F}_{n-1}(\Delta_t)}
      l^*(\Gamma),
    \end{split}
  \end{equation}
  where ${\mathcal F}_{n-1}(\Delta_t)$ is the set of
  $(n-1)$-dimensional faces of $\Delta_t$ and $l^*(\pi)$ denotes the
  number of lattice points that lie in the relative interior of the
  polytope $\pi$.
  
  The number of lattice points in a lattice polytope
  $\pi_t:=t\pi\subset \R^n$ is $E_\pi(t):= l(\pi_t)$, where $l(\pi_t)$
  is the number of lattice points in $\pi_t$ and $E_\pi$ is the
  Ehrhart polynomial of $\pi$ (see~\cite{ar:ehrhart-polynomials}).
  $E_\pi$ is a polynomial of degree at most $n$ where the coefficient
  of $t^n$ is $\vol(\pi)$. Furthermore, the reciprocity law says that
  $E_\pi(-t) = (-1)^n l^*(\pi_t)$.

  Since $D$ is an ample Cartier divisor, $\Delta$ is an
  $n$-dimensional lattice polytope~\cite[\S 2.2]{bo:oda-convex} and we
  have $l^*(\Delta_t) = (-1)^n E_\Delta(-t) = \vol(\Delta) t^n +
  O(t^{n-1})$ with $\vol(\Delta)>0$. Hence
  \begin{equation*}
    h^{n-1,0}_t \ =\ \vol(\Delta) t^n + O(t^{n-1}). 
  \end{equation*}
  Similarly,
  \begin{equation*}
    \begin{split}
      h^{n-2,1}_t \ &=\ l^*(2\Delta_t) -(n+1)l^*(\Delta_t) -
      \sum_{\Gamma\in {\mathcal F}_{n-1}(\Delta_t)} l^*(\Gamma) \\ &=\ 
      (-1)^n E_{2\Delta}(-t) - (n+1) (-1)^n E_{\Delta}(-t) -
      \sum_{\Gamma\in {\mathcal F}_{n-1}(\Delta_t)} (-1)^{n-1}
      E_\Gamma(-t) \\ &=\  \vol(2\Delta) t^n - (n+1) \vol(\Delta) t^n -
      \sum_{\Gamma\in {\mathcal F}_{n-1}(\Delta_t)} \vol(\Gamma)
      t^{n-1} + O(t^{n-1}) \\ &=\  (2^n - (n+1)) \vol(\Delta)t^n +
      O(t^{n-1}).
    \end{split}
  \end{equation*}
  
  Now we turn to $\mu_t$. By
  Theorem~\ref{th:toric_hypersurfaces_theory} we have $\mu_t = \dim
  R(f)_{[tD]}$ for any generic polynomial $f\in S_{[tD]}$. Since $\dim
  R(f)_{[tD]} = \dim S_{[tD]} - \dim J(f)_{[tD]}$, we study each term
  separately. Writing $D = \sum_{j=1}^r b_j D_j$ where the $D_j$ are
  the torus-invariant divisors associated to the one dimensional cones
  in the fan of $P$, we have
  \begin{equation*}
    \begin{split}
      \dim S_{[tD]} \ &=\ \# \{ (a_1,\ldots, a_r)\in \Z_{\geq 0}^r :
      \deg(z_1^{a_1}\cdots z_r^{a_r}) = [tD]\} \\ &=\ \# \{ m\in M :
      \langle m, n_j\rangle + tb_j \geq 0 \text{ for } j=1,\ldots, r\}
      \\ &=\ \# (\Delta_t \cap M) = l(\Delta_t),
    \end{split}
  \end{equation*}
  Also, since $J(f)_{[tD]} = \{ \sum_{j=1}^r g_j(z) \pd{f}{z_j} : g_j
  \in S_{[D_j]}, \forall j\}$, we have $\dim J(f)_{[tD]} \leq K$, for
  $K:= \sum_{j=1}^r \dim S_{[D_j]}$. Notice that $K$ is independent of
  $t$ (in fact, it is determined by the $1$-dimensional cones of the
  fan of $P$). We conclude then that
  \begin{equation*}
    \begin{split}
      \mu_t \ &=\ \dim R(f)_{[tD]} = \dim S_{[tD]} - \dim J(f)_{[tD]} \geq
      l(\Delta_t) -K \\ &=\ \vol(\Delta) t^n + O(t^{n-1}).
    \end{split}
  \end{equation*}
\end{proof}

\begin{proposition}\label{prop:dimensions_toric}
  Let $\beta\in A_{n-1}(P)$ be the class of an ample Cartier divisor
  in the complete simplicial toric variety $P$ of dimension $n\geq 4$.
  Then, there exists $t_0\in \N$ such that the cohomology of the
  generic ample hypersurface $X$ of degree $t\beta$ for $t\geq t_0$
  satisfies
  \begin{equation}\label{eq:comparison}
    \mu_t \ \geq\ 3 \left(\left[
    \frac{h^{n-2,1}_{v}(X) - 1}{h^{n-1,0}_{v}(X)}\right]  +1
    \right),   
  \end{equation}
  where $\mu_t:=\dim M_{t\beta}$.
\end{proposition}

\begin{proof}
  Let $D$ be an ample Cartier divisor with $[D]=\beta$. Then, by
  Lemma~\ref{le:dimensions_toric} we have 
  \begin{equation*}
    \begin{split}
      \left[ \frac{h^{n-2,1}_t - 1}{h^{n-1,0}_t}\right] \ &\leq\ 
      \frac{h^{n-2,1}_t - 1}{h^{n-1,0}_t} \ \leq\ 
      \frac{(2^n-(n+1))\vol(\Delta) t^n + O(t^{n-1})}{\vol(\Delta) t^n
        + O(t^{n-1})} \\ &=\ \frac{(2^n-(n+1))\vol(\Delta) +
        \frac{O(t^{n-1})}{t^n}}{\vol(\Delta) + \frac{O(t^{n-1})}{t^n}}
    \end{split}
  \end{equation*}
  and, since the last expression converges to $2^n-(n+1)$ as
  $t\rightarrow \infty$, we see that the whole expression is bounded
  from above.

  From
  \begin{equation*}
    H^{n-1}(X_f,\Q)_p \ =\ H^{n-1}(X_f,\Q)_v \oplus j^*H^{n-1}(P,\Q)
  \end{equation*}
  (\cite[Prop. 2.27]{bo:Voisin_Hodge_II}), Lefschetz's theorem and
  Bott's formula (see \cite[Thm.  2.14]{ar:materov-bott}), it is easy
  to conclude that, for $n\geq 4$, $H^{n-1,0}(X) = H^{n-1,0}(X)_v$ and
  $H^{n-2,1}(X) = H^{n-2,1}(X)_v$. Hence the Hodge numbers
  $h^{n-1,0}_t$ and $h^{n-2,1}_t$ compute the corresponding dimensions
  of the vanishing cohomology of $X$. 
  
  Thus, the right hand side of~\eqref{eq:comparison} is bounded from
  above. On the other hand, also from Lemma~\ref{le:dimensions_toric},
  we have $\mu_t \geq \vol(\Delta) t^n + O(t^{n-1})$ with
  $\vol(\Delta)>0$, so that the left hand side
  of~\eqref{eq:comparison} grows like $t^n$, and the inequality holds
  for $t$ sufficiently large.
\end{proof}

\begin{remark}
  In the proof of Proposition~\ref{prop:dimensions_toric} a choice of
  divisor $D$ is used. This choice determines how large $t$ should be
  so that~\eqref{eq:comparison} holds. But any other divisor $D'$ with
  the same degree would have an associated polytope $\Delta'$ that is
  a lattice translate of that of $D$; in particular, this implies that
  they have the same number of lattice points and hence the left and
  right sides of~\eqref{eq:comparison} would be the same for $D$ and
  $D'$. Hence, a value of $t$ that makes~\eqref{eq:comparison} true
  for $D$, works for any other $D'$ with degree $\beta$.
\end{remark}

Next we study the $p_1$ projection of the infinitesimal variation of
Hodge structure associated to sufficiently ample hypersurfaces in
toric varieties. We first recall, in a slightly adapted form, Lemma
1.28 of~\cite{ar:green-high_degree}:
\begin{lemma}\label{le:1.28}
  Let $P$ be a projective variety and ${\mathcal E}_1$, ${\mathcal
    E}_2$ two coherent sheaves over $P$. For ${\mathcal L}$ a
  sufficiently ample invertible sheaf, the multiplication map
  \begin{equation*}
    H^0(P,{\mathcal E}_1\otimes {\mathcal L}^a) \otimes 
    H^0(P,{\mathcal E}_2\otimes {\mathcal L}^b) \rightarrow 
    H^0(P,{\mathcal E}_1\otimes {\mathcal E}_2 \otimes {\mathcal L}^{a+b})
  \end{equation*}
  is surjective for $a,b\geq 1$.
\end{lemma}

\begin{remark}
  The ampleness condition in Lemma~\ref{le:1.28} is that 
  \begin{equation}
    \label{eq:condition-1.28}
    H^1(P\times P, {\mathcal I}_\Delta \otimes {\mathcal E}_1 \otimes 
    {\mathcal E}_2 \otimes \pi_1^* {\mathcal L}^a \otimes 
    \pi_2^* {\mathcal L}^b) \ =\ 0
  \end{equation}
  for all $a,b\geq 1$, where ${\mathcal I}_\Delta$ is the ideal sheaf
  of the diagonal in $P\times P$.
\end{remark}

\begin{proposition}\label{prop:p1=0}
  Let $P$ be a complete simplicial projective toric variety of
  dimension $n$ and $\beta\in A_{n-1}(P)$ an ample Cartier degree.
  Then there is a $t_0\in \N$ with the property that if $E$ is the
  infinitesimal variation of Hodge structure associated to a generic
  $f\in S_{t\beta}$ for $t\geq t_0$, then $p_1(E) \neq \{0\}$.
\end{proposition}
\begin{proof}
  Let ${\mathcal E}_1:={\mathcal O}_P$, ${\mathcal E}_2:=\Omega_P^n$,
  $a:=1$ and $1\leq b\leq n-1$. If $[D]=\beta$ for an ample Cartier
  divisor $D$, then ${\mathcal L}:={\mathcal O}_P(tD)$ will satisfy
  condition~\eqref{eq:condition-1.28} for all $t\geq t_0$ for some
  $t_0\in \N$. Using Lemma~\ref{le:1.28} with ${\mathcal E}_1$,
  ${\mathcal E}_2$, $a$ and $b$ as above and ${\mathcal L}:={\mathcal
    O}_P(tD)$ for $t\geq t_0$ we see that for a generic $f\in S_{t\beta}$
  \begin{equation*}
    H^0(P,{\mathcal O}_P(X_f)) \otimes H^0(P,\Omega_P^n(b X_f)) \rightarrow 
    H^0(P,\Omega_P^n((b+1) X_f)) 
  \end{equation*}
  is surjective. Then, by~\eqref{diag:compatibility_toric_2},
  \begin{equation}\label{eq:map}
    T_fB_\beta \otimes H^{n-p,p-1}(X_f)_v \rightarrow H^{n-p-1,p}(X_f)_v
  \end{equation}
  is surjective. Therefore, the map obtained by iterating $n-1$
  times~\eqref{eq:map}
  \begin{equation*}
    (T_fB_\beta)^{\otimes (n-1)} \otimes H^{n-1,0}(X_f)_v 
    \rightarrow H^{0,n-1}(X_f)_v
  \end{equation*}
  is surjective. But, since by taking $t_0$ larger if necessary
  $h^{n-1,0}\neq 0$ by Lemma~\ref{le:dimensions_toric}, we obtain that
  $p_1(E)\neq 0$.
\end{proof}

\begin{theorem}\label{th:non_genericity_toric}
  Let $P$ be a complete simplicial projective toric variety of
  dimension $n\geq 4$ and $\beta\in A_{n-1}(P)$ the degree of an ample
  Cartier divisor. For $t\in \N$ let $E_{t\beta}$ be the infinitesimal
  variation of Hodge structure on $M_{t\beta}$
  (see~\eqref{diag:factorization_through_moduli}) associated to a
  generic $f\in S_{t\beta}$. Then, there exists a $t_0\in \N$ such
  that, if $E_{t\beta}$ satisfies
  \begin{equation}\label{eq:conditions-nongeneric_toric}
    \dim(p_0(E_{t\beta})) \ =\ \dim E_{t\beta}, 
  \end{equation}
  for $t\geq t_0$, $E_{t\beta}$ is non-generic in $V_{\dim
    E_{t\beta}}$.
\end{theorem}

\begin{proof}
  This is immediate from Propositions~\ref{prop:dimensions_toric}
  and~\ref{prop:p1=0}, and Theorem~\ref{th:non-genericity-principle}.
\end{proof}

\begin{example}\label{ex:hypersurfaces_in_Pn}
  If $P:=\CP^n$, the infinitesimal Torelli and Macaulay's theorems
  imply that condition \eqref{eq:conditions-nongeneric_toric} holds.
  Fixing $D$ as the hyperplane divisor, we conclude that the
  infinitesimal variation
  associated to hypersurfaces of sufficiently high degree are non
  generic, that is, we recover the result
  of~\cite{ar:allaud-hypersurfaces_dn}.
\end{example}

We can now extend the result of the previous example to the case of
variations of generic smooth hypersurfaces in smooth toric varieties.  To
show that in this case~\eqref{eq:conditions-nongeneric_toric} holds we
quote Theorem 0.1 from~\cite{ar:green-high_degree}:

\begin{theorem}\label{th:0.1}
  Let $P$ be a smooth complete algebraic variety of dimension $n\geq
  2$ and $X$ a smooth member of the linear system determined by a
  sufficiently ample line bundle on $P$. Then, the map
  \begin{equation*}
    H^1(X,{\mathcal T}_X) \rightarrow \hom(H^0(X,\Omega^{n-1}_X), 
    H^1(X,\Omega^{n-2}_X))
  \end{equation*}
  is injective.
\end{theorem}

Then, in the smooth case, we conclude from
Theorems~\ref{th:non_genericity_toric} and~\ref{th:0.1} that

\begin{theorem}\label{th:smooth_main}
  Let $P$ be a complete smooth projective toric variety of dimension
  $n\geq 4$ and $\beta\in A_{n-1}(P)$ the degree of an ample Cartier
  divisor. Let $E$ be the infinitesimal variation of Hodge structure
  associated to a generic $f \in S_{t\beta}$. Then, there exists
  $t_0\in \N$ such that for all $t\geq t_0$, $E$ is non-generic.
\end{theorem}
\begin{proof}
  Let $t_0$ be the value produced by
  Theorem~\ref{th:non_genericity_toric}, and $X_f$ the hypersurface
  associated to a generic $f\in S_{t\beta}$ for $t\geq t_0$.
  Generically, $X_f$ is smooth. Since by increasing $t_0$ if needed,
  $T_f M_{t\beta} \hookrightarrow H^1(X_f,{\mathcal T}_{X_f})$, by
  Theorem~\ref{th:0.1}, we conclude that $E\simeq p_0(E)$ so
  that~\eqref{eq:conditions-nongeneric_toric} holds and so does this
  result.
\end{proof}

\begin{remark}
  Assuming an infinitesimal Torelli theorem valid for (generic)
  hypersurfaces of sufficiently high degree in toric varieties, then
  it is possible to remove the smoothness requirement from
  Theorem~\ref{th:smooth_main}. In the next section we show that
  smoothness is not an essential ingredient for the non-genericity
  result.
\end{remark}


\section{Hypersurfaces in Weighted Projective Spaces}
\label{sec:hypersurfaces_in_wpn}

Let $(q_0,\ldots,q_n)\in \N^{n+1}$ such that $q_0 = \gcd(q_1,\ldots,
q_n) = 1$ and $m|s$ where $m:=\lcm(q_1,\ldots,q_n)$ and
$s:=\sum_{j=0}^n q_j$.

Let $P:=\CP^n(q_0,\ldots,q_n)$ be the weighted projective space with
weights $(q_0,\ldots,q_n)$ (so that $P$ is ``well formed'' by the
conditions on the weights). $P$ is a complete simplicial toric variety
of dimension $n$ with fan $\Sigma := \{ \sigma_I : I \subset
\{0,\ldots,n\}\}$, where $\sigma_{\{i_1,\ldots,i_k\}}$ is the cone
generated by $e_{i_1}, \ldots, e_{i_k}$, with $\{e_1,\ldots,e_n\}$ the
canonical basis of $\R^n$ and $e_0 := -\sum_{j=1}^n e_j$. In
particular, the $1$-dimensional cones are $\Sigma^{(1)} :=
\{\sigma_{j} : j=0,\ldots n\}$ where $\sigma_j:=\R_{\geq 0} e_j$ and
the $n$-dimensional cones are $\Sigma^{(n)} := \{ \widehat{\sigma_j} :
j=0,\ldots n\}$, where $\widehat{\sigma_j}$ is the cone generated by
all the $e_l$'s except for $e_j$.  The lattice $N\subset \R^n$ is
generated by $\{f_1:=\frac{1}{q_1}e_1, \ldots,
f_n:=\frac{1}{q_n}e_n\}$, and we define $f_0:=e_0 = -\sum_{j=1}^n q_j
f_j$. Notice that $f_j$ is the integral generator of $\sigma_j$ for
$j=0,\ldots,n$.

$A_{n-1}(P) \simeq \Z$ and if $D_j$ is the torus-invariant divisor
associated to the cone $\sigma_j$, $\deg(D_j) = q_j$. Hence the
grading on $S:=\C[z_0,\ldots,z_n]$ is given by
\begin{equation*}
  \deg(z_0^{a_0} \cdots z_n^{a_n}) \ :=\ \sum_{j=0}^n q_j a_j  
\end{equation*}

\begin{proposition}
  Let $D:=\sum_{j=0}^n a_j D_j$ be a Weil divisor in $P$ with
  $d:=\deg(D)$. Then, $D$ is a Cartier divisor if and only if $m|d$.
  Also, $D$ is ample if and only if $d>0$.
\end{proposition}

\begin{proof}
  It is known (see~\cite{bo:fulton-toric}, Exercise on page 62) that
  $\sum_{j=0}^n a_j D_j$ is Cartier if and only if for each maximal
  cone $\sigma$ there is $u(\sigma)\in M$ such that for all $f_j$,
  lattice generators of the $1$-dimensional cones contained in
  $\sigma$,
  \begin{equation}
    \label{eq:property_of_u}
    \langle u(\sigma),f_j\rangle \ =\ -a_j.
  \end{equation}
  If $D$ is Cartier, there are $u(\widehat{\sigma_k})\in M$ for each
  $n$-dimensional cone $\widehat{\sigma_k}$, $k=0,\ldots, n$ that
  satisfy~\eqref{eq:property_of_u}. Evaluating each
  $u(\widehat{\sigma_k})$ on the lattice generators $f_j$ and
  using~\eqref{eq:property_of_u} it is easy to see that
  \begin{equation}
    \label{eq:definition_of_u}
    u(\widehat{\sigma_k}) \ =\ -\left(\sum_{j=1}^n a_j f_j^*\right) + \alpha_k 
    \quad \text{ with } \quad \alpha_k \ :=\
    \begin{cases}
      0 \text{ if } k=0\\
      \frac{d}{q_k} f_k^* \text{ otherwise}.
    \end{cases}
  \end{equation}
  Hence, since $u(\widehat{\sigma_k})\in M$, $q_k|d$ for all $k$, and
  so $m|d$.

  Conversely, if $m|d$, for each maximal cone
  $\widehat{\sigma_k}$~\eqref{eq:definition_of_u} defines an element
  of $M$ that satisfies~\eqref{eq:property_of_u}. Therefore, $D$ is
  Cartier.

  By Theorem 4.6 in~\cite{ar:batyrev-quantum}, $D$ is ample if and
  only if its piecewise linear support function $\psi_D$ satisfies
  \begin{equation}
    \label{eq:ampleness_condition}
    \psi_D(f_{i_1} + \cdots + f_{i_k}) - (\psi_D(f_{i_1}) + \cdots +
    \psi_D(f_{i_k})) \ >\ 0
  \end{equation}
  for all primitive collections $\{f_{i_1},\ldots,f_{i_k}\}$. In our
  case, the only primitive collection is $\{f_{0},\ldots,f_{n}\}$, for
  which the left hand side of~\eqref{eq:ampleness_condition} evaluates
  to $\frac{d}{\max \{q_0,\ldots,q_n\}}$ and the second assertion
  follows.
\end{proof}

\begin{corollary}
  $m D_0$ is an ample Cartier divisor. In particular, $P$ is a
  projective variety.
\end{corollary}

Let $E$ be the infinitesimal variation of Hodge structure associated
to a generic hypersurface of degree $d=tm$ for some $t\in \N$.  We are
interested in the dimension of $p_0(E)$, which will be computed using
the weighted Macaulay's theorem. Let us state it here, in a slightly
simplified form. Let $\rho:=(n+1)d-2s$,

\begin{theorem}[Weighted Macaulay's Theorem \cite{ar:tu-macaulay}]
  \label{th:macaulay} 
  If $c$ is a multiple of $m$ and $e \in \N$ such
  that $\rho-(c+e)>-s+mn$ then the natural map $R_e \rightarrow
  \hom(R_c,R_{c+e})$ induced by the ring multiplication is injective.
\end{theorem}

We can then derive the following

\begin{corollary} \label{cor:corollary:projections_wh} 
  If $d \geq 2m$ we have $\dim p_0(E)=\dim E$.
\end{corollary}

\begin{proof}
  To determine $\dim p_0(E)$ we apply Theorem~\ref{th:macaulay} for
  $e=d$ and $c=d-s$. Then our assumptions ensure that $m$ divides $c$.
  Moreover
  \begin{equation*}
    \rho-(c+e)\ =\ (n+1)d-2s-d-d+s\ =\ (n-1)d-s  
  \end{equation*}
  then as $d \geq 2m$ we have that $\rho-(c+e) \geq 2m(n-1)-s \geq
  mn-s$. The weighted Macaulay theorem then gives us that $R_d
  \hookrightarrow \hom(R_{d-s},R_{2d-s})$. Using the
  diagram~\eqref{diag:compatibility_toric} this gives us the following
  injection:
  \begin{equation*}
    E \hookrightarrow \hom(H^{n-1,0}, H^{n-2,2})  
  \end{equation*}
  which exactly means that $p_0(E) \simeq E$.
\end{proof}

\begin{theorem}\label{th:weighted_main}
  Let $P=\CP^n(q_0,\ldots,q_n)$ with $n\geq 4$, $q_0 =
  \gcd(q_1,\ldots, q_n) = 1$. Then, if $m|d$ and $m|s$, there is a
  degree $d_0$ such that for all $d\geq d_0$ the infinitesimal
  variation of Hodge structure associated to a generic hypersurface of
  degree $d$ is non-generic.
\end{theorem}

\begin{proof}
  Recall that $E$ is the infinitesimal variation of Hodge structure
  associated to a generic hypersurface of degree $tm$. By
  Corollary~\ref{cor:corollary:projections_wh}, if $t\geq 2$
  condition~\eqref{eq:conditions-nongeneric_toric} holds. On the other
  hand, since $m$ is the degree of an ample Cartier divisor, by
  Theorem~\ref{th:non_genericity_toric} the result follows.
\end{proof}



\section{Complete intersections}
\label{sec:ci_case}

In this section we prove that the infinitesimal variations of Hodge
structure associated to generic complete intersections of sufficiently
high degree in $\CP^n$ are non-generic. The context is the following:
we have $X\subset \CP^n$ given as the complete intersection of $c$
hypersurfaces of degrees $d_1\geq \cdots \geq d_c$, each determined by
a homogeneous polynomial $F_a$ for $a=1,\ldots, c$. We will assume
that $X$ is smooth. The degree of the canonical bundle of $X$, $K_X$,
will be denoted by $d(X)$, so that, $d(X):=\sum_{a=1}^c d_a -(n+1)$.
In this section we use the \emph{primitive} Hodge numbers, that is,
$h^{p,q} := \dim H^{p,q}(X)_p$.

The following Theorem states the relevant known results on the
cohomology of a complete intersection.

\begin{theorem}\label{th:ci_properties}
  With the notation as above, let
  $S:=\C[x_0,\ldots,x_n,\mu_1,\ldots,\mu_c]$, $F:=\sum_{a=1}^c \mu_a
  F_a \in S$, $J:=\langle \pd{F}{\mu_1}, \ldots, \pd{F}{\mu_c},
  \pd{F}{x_0}, \ldots, \pd{F}{x_n} \rangle$ the Jacobian ideal
  of $F$ and $R:=S/J$ seen as a bigraded ring with the grading induced
  by $\deg(x_j):=(0,1)$ and $\deg(\mu_a):=(1,-d_a)$.  Then:
  \begin{enumerate}
  \item \label{th:prop_cohomology} $H^{n-c-p,p}(X)_p \simeq
    R_{(p,d(X))}$, where $H_p$ is the primitive cohomology. In
    particular, $H^{n-c,0}(X)_p \simeq R_{(0,d(X))}$ and
    $H^{n-c-1,1}(X)_p \simeq R_{(1,d(X))}$.
  \item \label{th:prop_moduli} $H^1(X,{\mathcal T}_X) \simeq
    R_{(1,0)}$, where ${\mathcal T}_X$ is the tangent sheaf of $X$.
    Also, if $d(X)\geq 0$, the tangent space to the moduli of $X$ is
    isomorphic to $H^1(X,{\mathcal T}_X)$.
  \item \label{th:prop_compatible} The diagram
    \begin{equation} \label{eq:diag_compatibility}
      \xymatrix{{H^1(X,{\mathcal T}_X)} \ar[r] \ar[d]_{\simeq} &
      {\hom\left(H^{n-c-p,p}(X)_p,H^{n-c-p-1,p+1}(X)_p\right)} \ar[d]^{\simeq}
      \\ {R_{(1,0)}} \ar[r] & {\hom\left(R_{(p,d(X))},R_{(p+1,d(X))}\right)}}
    \end{equation}
  \end{enumerate}
  ---where the top horizontal arrow is the action of the infinitesimal
  variation and the bottom one is induced by multiplication---
  commutes up to a multiplicative constant.
\end{theorem}

Part~\eqref{th:prop_cohomology} is due to Terasoma, Konno,
Dimca~\cite[\S 2]{ar:dimca-complete_intersections}.  The first part
of~\eqref{th:prop_moduli} and part~\eqref{th:prop_compatible} are obtained
by adapting Terasoma's proofs of Propositions 2.5 and 2.6
in~\cite{ar:terasoma_weak_torelli} to the bigraded context. These
results are proved using a \emph{Cayley trick}, that is, by associating
to $X$ a hypersurface ${\mathcal X} :=(F=0)$ in a larger ambient space
and studying the relation between the cohomologies of $X$ and
${\mathcal X}$. The second assertion of~\eqref{th:prop_moduli} is
Theorem 3.5 in~\cite{ar:peters_l_torelli_I}.

\begin{remark}
  The complete intersection assumption ensures that the sequence
  $(F_1,\ldots, F_c)$ is regular.
\end{remark}

In order to apply Theorem~\ref{th:non-genericity-principle}, we need
to meet condition~\eqref{eq:ineq-principle} and check that the
conditions on the projections of the infinitesimal variation hold.
The purpose of the following result is to solve the first issue.

\begin{proposition}\label{proposition:bound_complete}
  Assume that $\dim X = n-c\geq 2$. Then there is $K_1\in \N$ such
  that, for $d(X)\geq K_1$, 
  \begin{equation}\label{eq:inequality}
    \mu\ \geq\ 3\left( \left[ \frac{h^{n-c-1,1} -1}{h^{n-c,0}} \right] +
      1\right), 
  \end{equation}
  where $\mu$ is the dimension of the infinitesimal variation of Hodge
  structure associated to the family of deformations of $X$. 
\end{proposition}

\begin{proof}
  This proposition follows readily from polynomial estimates of the
  different dimensions that appear in~\eqref{eq:inequality}.

  Using combinatorics, it is not difficult to show the following
  bounds:
  \begin{equation*}
    \begin{split}
      h^{n-c-1,1} \ &\leq\ \frac{c}{n!} d(X)^n+  O(d(X)^{n-1})\\
      \dim H^1(X,{\mathcal T}_X) \ &\geq\ \frac{1}{n!  c^n} d(X)^n + O(1)
    \end{split}
  \end{equation*}
  and, for $c\geq 2$, 
  \begin{equation*}
    h^{n-c,0}\ \geq\ \frac{1}{n^n} d(X)^{n-c}
  \end{equation*}
  while, for $c=1$,
  \begin{equation*}
    h^{n-1,0} \ \geq\ \frac{1}{n!} d(X)^n.
  \end{equation*}
  For $c\geq 2$, the right hand side of~\eqref{eq:inequality} is
  bounded from above by $3 \frac{n^n c 2^n}{n!}  d(X)^c +
  O(d(x)^{c-1})$, and, $\mu = \dim H^1(X,{\mathcal T}_X) \geq
  \frac{1}{n!  c^n} d(X)^n + O(1)$. Then, since $c<n$, the statement
  follows. The case for $c=1$ is the analogous.
\end{proof}

\begin{remark}
  It is possible to give an effective version of
  Proposition~\ref{proposition:bound_complete}. Indeed, if   
  \begin{equation*}
    d(x)\ \geq\ 
    \begin{cases}
      \max\{n,\sqrt[n-c]{3n^n c^{n+1} 2^n +1},\sqrt[c]{n!c^n
        (3+c^2+(n+1)^2)} \}, \text{ if } c\geq 2,\\
      \max\{ n, \sqrt[n]{n!(2^n 3 + 4+(n+1)^2)}\}, \text{ if } c=1
    \end{cases}
  \end{equation*}
  then~\eqref{eq:inequality} and, Theorem~\ref{th:ci_main} hold. Since
  we don't find this explicit bound very illuminating, or useful, we
  omit the details.
\end{remark}

Next we prove that the conditions on the projections of the
infinitesimal variation required by
Theorem~\ref{th:non-genericity-principle} hold.

\begin{proposition} \label{proposition:projections} There is a
  $K_2\in\N$ such that if $E$ is a complete intersection \IVHS\ with
  $d(X)\geq K_2$, then we have:
  \begin{equation*}
    p_0(E) \ \simeq\ E \quad \text{ and }\quad  p_1(E) \ \neq\ \{0\}    
  \end{equation*}
\end{proposition}

\begin{proof}
  First we refer to the local Torelli theorem proved by C. Peters
  in~\cite{ar:peters_l_torelli_I}. There he proved that the following
  map is injective:
  \begin{align*}
    E &\hookrightarrow \hom(H^{n-c,0}(X)_p,H^{n-c-1,1}(X)_p)\\
    \alpha &\longmapsto \alpha_{\mid H^{n-c,0}(X)_p}
  \end{align*}
  This means exactly that $E \simeq p_0(E)$.  Next, using Lemma 3.4 of
  J. Nagel~\cite{ar:nagel-complete_intersections}, we see that the map
  \begin{equation*}
    R_{(1,0)} \otimes R_{(1,d(X))} \longrightarrow R_{(2,d(X))}    
  \end{equation*}
  induced by multiplication is surjective, so using
  diagram~\eqref{eq:diag_compatibility} we find that the corresponding
  map in cohomology $E \otimes H^{n-c-1,1}(X)_p \rightarrow
  H^{n-c-2,2}(X)_p$ is surjective; therefore $p_1(E) \neq \{0\}$ as
  soon as $h^{n-c-2,2}$ is non-zero. But using the surjectivity of the
  multiplication in $R$ and diagram~\eqref{eq:diag_compatibility} once
  more, the following map is also surjective:
  \begin{equation*}
      E \otimes H^{2,n-c-2}(X)_p \longrightarrow H^{1,n-c-1}(X)_p.
  \end{equation*}
  Now, as seen in the proof of
  Proposition~\ref{proposition:bound_complete}, $\dim H^1(X,{\mathcal
    T}_X) \geq \frac{1}{n!  c^n} d(X)^n + O(1)$, so that there is
  $K_2\in\N$ with the property that $H^1(X,{\mathcal T}_X) \neq 0$ if
  $d(X)\geq K_1$. If this last condition holds, then $E\neq \{0\}$.
  Therefore, $h^{1,n-c-1}=h^{n-c-1,1}$ is non-zero because $E$ injects
  into $\hom(H^{n-c,0}(X)_p,H^{n-c-1,1}(X)_p)$; this then implies that
  $h^{n-c-2,2}=h^{2,n-c-2}\neq 0$.
\end{proof}

From Theorem~\ref{th:non-genericity-principle}, whose hypotheses are
satisfied because of Propositions~\ref{proposition:projections}
and~\ref{proposition:bound_complete}, we derive the final result:

\begin{theorem}\label{th:ci_main}
  Let $X\subset \CP^n$ be a smooth complete intersection defined by
  the intersection of $c$ hypersurfaces of degrees $d_1,\ldots,d_c$.
  If $\dim X \geq 3$ there exists $K\in\N$ such that if the degree of
  the canonical bundle $d(X) = \sum_{a=1}^c d_a -(n+1)\geq K$, then
  the \IVHS associated to the deformation of $X$, $E$, is non-generic
  in $V_{\dim E}$. 
\end{theorem}





\def\cprime{$'$}
\providecommand{\bysame}{\leavevmode\hbox to3em{\hrulefill}\thinspace}
\providecommand{\MR}{\relax\ifhmode\unskip\space\fi MR }
\providecommand{\MRhref}[2]{%
  \href{http://www.ams.org/mathscinet-getitem?mr=#1}{#2}
}
\providecommand{\href}[2]{#2}


\end{document}